\input amstex
\documentstyle{amsppt}
\mag=\magstep1
\document

\topmatter
\title
Truncated Counting Functions of Holomorphic Curves 
in Abelian Varieties
\endtitle
\author
Ryoichi Kobayashi
\endauthor
\affil
Graduate School of Mathematics, Nagoya University
\endaffil
\abstract{A new proof of the Second Main Theorem with 
truncation level $1$ for Zariski-dense holomorphic curves into Abelian 
varieties, which has just been proved by Yamanoi [Y2], is presented. 
Our proof is based on the idea of the ``Radon transform" introduced in [K2] 
combined with consideration on certain singular perturbation of the probability 
measures on the parameter space which appears in the ``Radon transform".}
\endabstract
\endtopmatter
\rightheadtext{Truncated Counting Functions}

\def\C{\Bbb C} \def\Euc{\text{\rm Euc}} \def\Ram{\text{\rm Ram}} \def\mult{\text{\rm mult}}
\def\Supp{\text{\rm Supp}}   
\def\ep{\varepsilon} \def\NI{\newline\noindent} \def\PP{\Bbb P} \def\l{\lambda}
\def\L{\Lambda} \def\Hilb{\text{\rm Hilb}} \def\mer{\cdots\to} \def\om{\omega} \def\wt{\widetilde}
\def\Q{\Bbb Q} \def\Z{\Bbb Z} \def\p{\partial} \def\dist{\text{\rm dist}} \def\exc{/\kern-3pt/\kern-1pt}
\def\Sing{\text{\rm Sing}}  \def\p{\partial}

\beginsection 1. Introduction.

The Nevanlinna Theory describes ``intersection theory" of 
holomorphic curves from $\C$ and divisors in a projective 
algebraic variety in terms of the relationship among ``basic 
functions". Let us recall the definition of the basic functions 
in Nevanlinna theory. Let $X$ be a smooth projective 
variety, $D$ an effective divisor, $\sigma$ a canonical section 
of $\Cal O_X(D)$ and $\Vert\cdot\Vert$ a smooth Hermitian norm 
of $\Cal O_X(D)$. Let $f:\C \to X$ be a holomorphic 
curve such that $f(\C) \not\subset \Supp(D)$. 
The ``intersection theory" described by the Nevanlinna theory 
of (transcendental) holomorphic curves and divisors 
has two aspects. One is to measure how a holomorphic curve can 
approximate a given divisor. The other is to measure how often a 
holomorphic curve can intersect a divisor. 
We measure the approximation of a holomorphic curve $f:\C\to X$ by the asymptotic 
behavior of the {\bf proximity function}
$$
m_{f,D}(r)=\int_0^{2\pi}\log^+\frac1{\Vert\sigma\circ f(re^{i\theta})\Vert}
\frac{d\theta}{2\pi}\,\,.
$$
Note that $\Vert\sigma\circ f(re^{i\theta})\Vert$ is equivalent to the Euclidean distance 
between $f(re^{i\theta})$ and $D$. 
Let $n_{f,D}(t)$ (resp. $n_{f,D}(0)$) denote the numbers of zeros of $\sigma\circ f$ 
in $\C(t)=\{z\in\C\,;\,|z|<t\}$ (resp. $\deg_0(\sigma\circ f)$). 
We then measure how often $f$ intersects $D$ by the asymptotic behavior of 
the {\bf counting function}
$$
\split
N_{f,D}(r) & =\sum_{0\not=a\in\C(r)}\deg_a(\sigma\circ f)
\log\biggl|\frac{r}{a}\biggr|+n_{f,D}(0)\log r\\
& = \int_0^r\frac{n_{f,D}(t)-n_{f,D}(0)}{t}dt
+n_{f,D}(0)\log r\,\,.
\endsplit
$$
The sum of these two functions is the {\bf height function}
$$
T_{f,D}(r)=m_{f,D}(r)+N_{f,D}(r)\,\,.
$$
We call these three functions the ``basic functions" in Nevanlinna theory. 

The First Main Theorem in Nevanlinna Theory states that the asymptotic 
behavior of the height function, considered as an element of 
$$\frac{\text{\rm all functions on $\Bbb R_{>0}$}}
{\text{\rm all bounded functions on $\Bbb R_{>0}$}}$$
depends only on the linear equivalence class of the divisor $D$. 
Indeed, if $D_1=(\sigma_1) \sim D_2=(\sigma_2)$ 
(linearly equivalent), then the difference of $T_{f,D_1}(r)-T_{f,D_2}(r)$ 
essentially depends on the leading coefficient of the Laurent expansion 
of the meromorphic function $\frac{\sigma_1}{\sigma_2}\circ f$ at $z=0$. 

There is another kind of counting function, namely, the {\bf ramification counting 
function}. Let $R$ be a Riemann surface and $f:\C \to R$ a holomorphic map. 
Then the ramification counting function is defined as
$$
N_{f,\Ram}(r)
=\sum_{0\not=a\in\C(r)}(\mult_a(f)-1)\log\biggl|\frac{r}{a}\biggr|
+(\mult_0(f)-1)\log r
$$
where $\mult_a(f)$ is the multiplicity of $f$ at $z=a$ (i.e., 1 plus 
the vanishing order of the derivative $f'$ at $z=a$). 
This coincides with the usual counting function of the solutions of the equation 
$f'(z)=0$ counted with multiplicities. 
In this article, the ramification counting function plays an essential role. 

The estimate of the Second Main Theorem type is most important 
in Nevanlinna Theory. For instance, the Conjecture 1.1 is a typical example 
of the leading conjecture in the Nevanlinna Theory 
(see, for instance, [L1,2], [NoO], [V1]):

\proclaim{Conjecture 1.1} Let $X$ be a smooth projective variety, $D$ a divisor 
with at worst simple normal crossings and $E \to X$ any ample line bundle. 
Let $\ep$ be any positive number. 
Then there exists a proper algebraic subset $Z=Z(X,D,\ep,E)$ with the following 
property. Let $f:\C \to X$ be any holomorphic curve such 
that $f(\C) \not\subset \Supp(D)$. Then if $f$ is algebraically non-degenerate in 
the sense that the the image of $f$ is Zariski dense, we have
$$
m_{f,D}(r)+N_{W(f),0}(r)+T_{f,K_X}(r) 
\leq \ep\,T_{f,E}(r)\exc_{\ep}\,\,,
$$
where $N_{W(f),0}(r)$ is the counting function for the equation $W(f)=0$, 
where $W(f)$ is a conjectural object interpreted as a kind of ``Wronskian'' 
determinant of the holomorphic curve $f:\C \to X$. Moreover, if $f$ is algebraically 
non-degenerate in the weaker sense that $f(\C) \not\subset Z$ holds, 
then there is a modification of $W(f)$ so that the inequality of the same form 
with modified $W(f)$ hold.
\endproclaim

The term $\ep\,T_{f,E}(r)$ is an error term and $\exc_{\ep}$ means that 
the inequality holds outside an $\ep$-dependent Borel set of $\ep$-dependent 
finite Lebesgue measure. 
In general, the counting function $N_{W(f),0}(r)$ should generalize 
the ramification counting function for holomorphic maps to Riemann surfaces. 

The following list essentially exhausts all known cases where Conjecture 1.1 
is true. 
\NI
(i) [Nevanlinna, Ahlfors] The case $X$ is a compact Riemann surface ([N], [A]). 
In this case, $Z=\emptyset$. 
\NI
(ii) [Nevanlinna-Cartan theory and its generalization] The case $X=\PP^{n}(\C)$, 
$D$ a collection of hyperplanes in general position and the holomorphic curve 
$f:\C \to \PP^n(\C)$ is linearly non-degenerate.  
([C], [W] , [A], [F] and [V2,3]). In [V2,3] Vojta proved that there exists a (effectively computable) 
finite collection of proper linear sybspaces $Z$ such that if $f(C)$ is not contained in $Z$, 
the inequality of Conjecture 1.1 holds for linearly non-degenerate $f$ and even for linearly 
degenerate $f$ such that $f(\C) \not\subset Z$ the same inequality without the Wronskian 
term (i.e., the term like $N_{W(f),0}(r)$) holds. It would be an interesting problem 
to recover a ``Wronskian term" in the linearly degenerate case. Recently, Ru [R] 
proved a defect relation when $D$ is a general hypersurface configuration. 
\NI
(iii) [Bloch, Ochiai] The case that $X$ is a subvariety (of general type) of an Abelian variety 
([O], [GG], [NoO], [KS], [Y2]). In this case, $D$ is empty and $Z$ is the union of all 
translations of proper Abelian subvarieties contained in $X$.
\NI
(iv) The case that $X$ is an Abelian variety and $D$ is an arbitrary divisor 
([SY], [Y1,2], [NoWY], [M], [K2]). In this case $Z=\emptyset$. 

These results are called the Second Main Theorem. The most well-understood 
Second Main Theorem for holomorphic curves into higher dimensional targets is 
the case (ii) where $X=\PP^{n}(\C)$ and $D$ a collection of $q$ linear divisors 
in general position. In this case $W(f)$ is the usual Wronskian determinant of $f$ 
with respect to the canonical affine local coordinate system of $\PP_n$. 
Replacing $\deg_a(\sigma\circ f)$ by $\min\{\deg_a(\sigma\circ f),k\}$ in the 
definition of the usual counting function, we get the counting function 
truncated at level $k$, which we denote by $N_{k,f,D}(r)$. Define the residual 
counting function $N^k_{f,D}(r)$ with truncation level $k$ by
$$
N^k_{f,D}(r)=N_{f,D}(r)-N_{k,f,D}(r)\,\,.
$$
So $N^k_{f,D}(r)$ counts only intersections of $f$ and $D$ with 
multiplicity $m \geq k$ with weight $m-k$ 
(i.e., replacing $\deg_a(\sigma\circ f)$ 
by $\max\{\deg_a(\sigma\circ f)-k,0\}$ in the definition of 
the usual counting function). We then have
$$
N^l_{f,D}(r) \leq N^k_{f,D}(r)
$$
if $k \leq l$. Restricting the equation $W(f)=0$ to $D$, we have
$$
N^n_{f,D}(r) = N_{W(f),0}(r)\,\,.
$$
Therefore, if $f:\C\rightarrow \PP^n(\C)$ is linearly non-degenerate, then 
the Second Main Theorem with the truncated counting function of truncation 
level $n$ holds:
$$
m_{f,D}(r)-(n+1)T_{f,\Cal O(1)}(r)+N^n_{f,D}(r) \leq O(\log^+(rT_{f,H}(r)))\exc\,\,.
$$
It follows from the above argument that the reason for the truncation 
level being $n=\dim\PP^{n}(\C)$ comes from the key role played by 
the Wronskian determinant $W(f)$. In general case, we conjecture that 
the Wronskian $W(f)$ must be replaced by its generalization (still some kind 
of Wronskian), and if so, the expected truncation level will be $n=\dim X$. 
\medskip

On the other hand, we expect some special property on the truncated counting 
function, if the target $X$ is a special variety, for instance, 
if $X$ is an Abelian variety. In fact, the theory of theta divisors of Jacobian varieties 
of compact Riemann surfaces suggests that the value distribution of 
holomorphic curves into the Jacobian variety of a compact Riemann surface 
should have some resemblance to the value distribution theory with one-dimensional 
targets. For instance, it is natural to ask if the relevant 
truncation level is one for holomorphic curves in Jacobian varieties (more generally, 
holomorphic curves in Abelian varieties). 
\medskip

Inspired by works of McQuillan [M] and Brunella [B], Yamanoi [Y2] 
discovered that the truncation level is taken to be one in the Second Main Theorem 
for {\it Zariski-dense} holomorphic curves $f:\C \to A$ into an Abelian variety $A$ ($D$ 
being any divisor):
$$
m_{f,D}(r)+N^1_{f,D}(r) \leq \ep\,T_{f,E}(r)\exc_{\ep}\,\,. \tag1
$$
The assumption $f$ being Zariski-dense cannot be removed. 
Moreover, the error term is of the form $\ep\,T_{f,E}(r)$ 
($\ep$ being any positive number) and this is not improved to 
the form $O(\log^+(rT_{f,E}(r)))$ ([Y2]). 
As an application of this result, Yamanoi [Y2] gave a new proof of 
the Bloch-Ochiai Theorem [O]. 
\medskip

The purpose of this article is to give a new proof of Yamanoi's 
result (1) ([Y2]) by using the method of ``Radon transform" developed in [K2].  
In the course of the proof, we will clarify a simple reason
\footnote{\,\,Definition 2.1-2.2 of the Radon transform and its property (response 
under certain singular perturbation of the probability space $\L$ involved in the 
Radon transform), which is presented in Lemmas 2.5 and 2.7.} 
why the truncation level in the Second Main Theorem is taken to be one 
for Zariski-dense holomorphic curves into Abelian varieties. 
\medskip

Our proof is separated into two steps. 

In the first step, we introduce the idea of the 
``Radon transform" of holomorphic curves in Abelian varieties ([K2]) in order to 
study the intersection of a given holomorphic curve $f:\C \to A$ with an ample 
divisor $D$. 
The ``Radon transform" transforms a given entire holomorphic 
curve $f:\C \to A$ into a family of holomorphic 
maps $\{f_{\l}:Y_{\l} \to S_{\l}\}_{\l\in\L}$, where $Y_{\l}$ is a finite analytic 
covering with projection $\pi_{\l}:Y_{\l} \to \C$ and $\{S_{\l}\}_{\l\in\L}$ form 
a family of algebraically equivalent algebraic curves in $A$ all passing through the 
neutral point. where $\L$ is the parameter space. 
To get information independent of $\l$, we apply the Second Main Theorem 
for holomorphic curves from finite covering of $\C$ to Riemann surfaces 
(see, for instance, [No]) to each $f_{\l}:Y_{\l}\rightarrow S_{\l}$ and ``average" 
them over the parameter space $\L$ with respect to some probability measure
\footnote{\,\,Which probability measure 
we should choose is essential in this article.}. 
In \S2, we introduce the ``Gaussian divisor" $D_{\l}^{[1]}\subset \PP(TA)$ (resp. the 
``incidence divisor" $D_{\l} \subset A$) which are associated to the pair $(D,S_{\l})$, 
which is related to the ramification of $f_{\l}:Y_{\l}\rightarrow S_{\l}$ 
(resp. $\pi_{\l}:Y_{\l}\rightarrow \C$). 
Then $D$, $D_{\l}$, $D_{\l}^{[1]}$ and $H^{[1]}$ (the relative hyperplane bundle 
of $\PP(TA)$) fits into the generalized Hurwitz formula (see (2) in \S2). 
We introduce a notion ``variation" of probability spaces $(\L,d\l_t)$ ($0<t<r$) with 
certain ``measure concentration", over which Nevanlinna theoretic 
functions containing $\l$ as a parameter are ``averaged". 
By examining the response of various counting functions under the 
``measure concentration", we know the special feature of ``averaging over $(\L,d\l_t)"$ 
procedure. 
For any $\ep>0$, there exists a variation of probability spaces $(\L,d\l_t)$ such that:
\smallskip
\noindent
$\bullet$ The ``average" of $N_{f_{\l},\Ram}(r)-N_{f^{[1]},D_{\l}^{[1]}}(r)$ is negligible.
\smallskip
\noindent
$\bullet$ The ``average" of $N_{\pi_{\l},\Ram}(r)-N_{f,D_{\l}}(r)$ is equal to $-N^1_{f,D}(r)$ 
modulo the small error of magnitude $\ep\,T_{f,D}(r)$. 
\smallskip

Combining these estimates with the generalized Hurwitz formula give rise to 
a Diophantine inequality (1) of the Second Main Theorem type, 
which is much stronger than the version with generally expected counting function 
$N^n_{f,D}(r)$ in the Second Main Theorem 
for holomorphic curves into $n$-dimensional Abelian varieties. 

The second step, where we still use special properties of Abelian targets, 
consists of showing that the proximity function of the jet lifts of a Zariski-dense 
holomorphic curve $f:\C \to A$ with respect to the jet spaces of $D$ 
is small in the sense that it is dominated by the quantity $\ep\,T_{f,E}(r)$ 
(for any $\ep>0$ and any ample line bundle $E\to A$). 
Nevanlinna's lemma on logarithmic derivative (see Lemma 3.1) applied to the higher 
jet lifts of $f:\C\to A$ is essential in the argument of the second part. 

The essential steps in our proof of (1) are Lemmas 2.5 and 2.7. 
To generalize these Lemmas to more general ``Radon transform" defined for 
arbitrary holomorphic curves into general projective algebraic varieties seems to 
be promissing because of the rich flexibility in the definition of the 
``Radon Transform" (see [K3]). We therefore pose the following conjecture for 
future study:

\proclaim{Conjecture 1.2} Let $X$, $D$, $E$ and $\ep$ be as in Conjecture 1.1. 
Let $f:\C \to X$ be any holomorphic curve which is algebraically non-degenerate 
in the sense that the the image of $f$ is Zariski dense in $X$. Then, we have
$$
m_{f,D}(r)+N^1_{f,D}(r)+N_{W(f),0}(r)+T_{f,K_X}(r) 
\leq \ep\,T_{f,E}(r)\exc_{\ep}\,\,,
$$
where $N_{W(f),0}(r)$ is a ``conjectural" counting function for the equation $W(f)=0$, 
where $W(f)$ is a ``conjectural" object interpreted as a kind of ``Wronskian'' 
determinant of the holomorphic curve $f:\C \to X$. Moreover, there exists a proper 
algebraic subset $Z=Z(X,D,\ep,E)$ with the following property. If $f:\C\to X$ is any 
holomorphic curve satisfying the property $f(\C) \not\subset Z$, then there is 
a modification of $W(f)$ so that the same inequality as above should hold.
\endproclaim

\beginsection 2. Radon Transformation.

Let $A$ be an Abelian variety and $D$ a reduced divisor. 
Suppose that $D$ is not ample. Then there exists an Abelian variety $A'$ 
of lower dimension, an ample divisor $D' \subset A'$ and a surjective 
morphism $\pi:A \to A'$ such that $\pi^*D'=D$. Therefore, the intersection theory 
of $f:\C\rightarrow A$ and $D$ is reduced to that of $\pi\circ f:\C\rightarrow A'$ 
with $D'$. By this reason, we put the 
following assumptions:
\smallskip
\noindent
$\bullet$ {\it Assumption 1.} The divisor $D$ is ample and reduced.
\smallskip

Moreover, we assume the following
\smallskip
\noindent
$\bullet$ The divisor $D$ has at worst normal crossings. 
\smallskip

This assumption is not necessary for the proof of (1) in the sense that our arguments 
can be modified to cover the cases without this. 
However, this assumption makes the arguments in this article considerably simple. 

Let $\{S_{\l}\}_{\l\in\L}$ ($\L$ being the parameter space) be an algebraic family 
of algebraically equivalent curves in $A$ all passing through the neutral point $p$ 
of $A$. 
\smallskip
\noindent
$\bullet$ {\it Assumption 2.} There is a Zariski open subset $\L^0$ of $\L$ 
such that  for all $\l\in\L^0$ the curve $S_{\l}$ is non-singular.
\smallskip
\noindent
$\bullet$ {\it Assumption 3.}  There is a positive integer $k$ with the following property. 
The natural rational map $t:\L \mer A^{[k]}_p$ 
from $\L$ to the $k$-th projective jet space $A^{[k]}_p$ at $p$ is surjective and 
generically finite. 

The $k$-th projective jet space $A^{[k]}_p$ at the neutral point $p\in A$ is defined 
as follows: we start with the $k$-times iterated projective tangent bundle of $A$, i.e., 
the space obtained from $A$ by iterating the operation of taking projective tangent 
bundle $k$-times. 
Then, $A^{[k]}_p$ consists of projective classes of $k$-jets of all germs of 
holomorphic curves passing through the neutral point $p$. 
The natural map $t$ is defined by
$$\L \ni \l \mapsto t(\l):=(S_{\l}^{[k]})_p \in A^{[k]}_p\,\,.$$ 
We introduce the measure 
$t^*(\text{\rm Fubini-Study measure of $A^{[k]}_p$})$ on $\L$ and call it the ``Fubini-Study 
measure". 
\smallskip
\noindent
$\bullet$ {\it Assumption 4.} Let $k$ be as in {\it Assumption 3} and consider 
the sequence $j=1,2,\dots,k$. 
Then, to each hyperplane in $T_pA$ 
is associated a strictly decreasing sequence of non-empty Zariski closed subsets 
$\{\L_j\}_{j=1}^k$ of $\L^0$ consisting of those $\l$ with the property that the local 
intersection number of $S_{\l}$ and $T_pA$ at $p$ is at least $j$.  
\medskip
\comment
From here on we assume that $D$ is a sufficiently ample divisor and we 
put the following assumption ($*$)\footnote{\,\,Choosing the family of curves $S_{\l}$ 
satisfying the assumption ($*$) is essential in the proof of Lemma 2.5.} : 
\medskip

($*$)\quad For each $a\in D$, the family $\L_a$ of curves $S_{\l}$ 
with the property that $(a+S_{\l}) \subset D$ forms a subvariety of $\L$ 
of codimension 1. 
It is almost equivalent to saying that if $a+S_{\l}$ is tangent to $D$ at $a$ 
then the whole $a+S_{\l}$ is contained in $D$.
\medskip
\endcomment

Later we will modify {\it Assumption 3 and 4} and the ``Fubini-Study measure" 
in the following way:
\medskip

1) We consider {\it Assumption 3 and 4} for various $k\in\Z_{>0}$. 
Namely, if we write $\L_k$ for such $\L$ satisfying the {\it Assumption 3 and 4} 
with positive integer $k$, we will consider families $\{S_{\l}\}_{\l\in\L_k}$ 
for various $k\in\Z_{>0}$. 
\medskip

2) We consider a certain variation of the measure 
$$t^*(\,\text{\rm ``Fubini-Study measure" of $[A]^{(k)}$})$$
on the parameter space $\L$. 
Namely, we will consider the variation of measures ``interpolating" 
the Fubini-Study measure and the Dirac measure ``perfectly concentrated" 
at some $\L_j$ associated to some hyperplane in $T_pA$. 
In particular, this variation of measures contains a regular measures which strongly 
``concentrate" at some $\L_j$ (w.r.to some hyperplane of $T_pA$). 
Such a variation of measures will play an essential role in our proof of (1). 
 \medskip

Let $k:=S_{\l}\cdot D$ (the intersection number) and $m=g(S_{\l})$ 
(the genus of $S_{\l}$, if $S_{\l}$ is smooth). For $a\in A$ and a subset $Z$ of $A$ 
we set
$$
Z+a=\{\zeta+a\,;\,\zeta\in Z\}\,\,.
$$
This is the translate of $Z$ by an element $a\in A$. 

\proclaim{Definition 2.1} Let $\phi_{\l}:A \to \Hilb^k(S_{\l})$ denote the 
holomorphic map 
$$\phi_{\l}(a)=(D-a)\cap S_{\l}$$
repeated according to multiplicities. The {\bf incidence divisor} 
associated to $D$ and $S_{\l}$ is defined by
$$
D_{\l}:=\phi_{\l}^*(\triangle)
$$
where $\triangle$ is the incidence divisor of $\Hilb^k(S_{\l})$.
\endproclaim

Note that $\Hilb^k(S_{\l})$ is a non-singular projective variety 
because $S_{\l}$ is a curve (for instance, the elementary symmetric polynomials 
induces the well-known isomorphism $\Hilb^k(\PP^1(\C))=\PP^k(\C)$). 
Let $f:\C \to A$ be an arbitrarily given holomorphic curve. 

\proclaim{Definition 2.2} (i) For each $\l\in\L$, we define
$$
Y_{\l}:=\text{\rm the\,\,normalization\,\,of\,\,}
\{(z,w) \in \C \times S_{\l}\,;\,w \in (D-f(z))\cap S_{\l}\}\,\,,
$$
and an analytic covering map
$$
\pi_{\l}\,:\,Y_{\l} \rightarrow \C\,\,;\qquad (z,w) \mapsto z
$$
of degree $k$ and a holomorphic map
$$
f_{\l}\,:\,Y_{\l} \rightarrow S_{\l}\,\,;\qquad (z,w) \mapsto w\,\,.
$$
\noindent
(ii) The {\bf Radon transformation} of a given holomorphic 
curve $f:\C \rightarrow A$ with respect to $\{S_{\l}\}_{\l\in\L}$ is defined 
to be the collection of holomorphic 
maps
$$\{f_{\l}:Y_{\l} \rightarrow S_{\l}\}_{\l\in\L}\,\,.$$
\endproclaim

The basic tool to count the ramification of $f_{\l}:Y_{\l}\ \rightarrow S_{\l}$ is the 
Gauss map $\C \ni z \mapsto [f'(z)] \in \PP(T_{f(z)}A)$. 
To identify the ramification counting function of $f_{\l}:Y_{\l}\to S_{\l}$ with a 
usual counting function with respect to a divisor, we introduce the notion of the 
Gaussian divisor:

\proclaim{Definition 2.2 (continued)} (i) The {\bf universal Gaussian divisor} 
$$D_{\L}^{[1]}
\subset (\PP(TA)\times \L)$$
associated to $D$ and the family $\{S_{\l}\}_{\l\in\L}$ is defined by
$$
\split
D_{\L}^{[1]} & := \text{\rm the\,\,Zariski\,\,closure\,\,in\,\,}
\PP(TA)\times\L\,\,\text{\rm of}\\
& \bigcup
\Sb (a,\l)\in A\times\L,\,\,D-a\,\,\text{\rm is\,\,non-singular}\\ 
\text{\rm at\,\,the\,\,intersection\,\,}(D-a)\cap S_{\l}\endSb 
\biggl(\bigcup_{w\in (D-a)\cap S_{\l}}
\underbrace{\PP(T_wD)}_{\text{\rm translated\,\,to\,\,}a}
\times\{\l\}\biggr)\,\,.
\endsplit
$$
\noindent
(ii) The {\bf Gaussian divisor} $D_{\l}^{[1]}$ associated to $S_{\l}$ is 
defined by
$$
D_{\l}^{[1]}:=D_{\L}^{[1]}\cap (\PP(TA)\times \{\l\})\,\,.
$$
\endproclaim

The Gaussian divisor $D_{\l}^{[1]}$ is a divisor of $\PP(TA)$ such that $D_{\l}^{[1]}
\cap \PP(T_aA)$ consists of $k$ hyperplanes in $\PP(T_aA)$ for general $a\in A$. 
These hyperplanes are the collection of the translates to $a$ of the projectivized 
tangent hyperplanes of $D-a$ at the intersection $(D-a)\cap S_{\l}$. 

Let $k:=S_{\l}\cdot D$. We have a $k$-valued holomorphic ``map" 
$F_{D,\l}:A\rightarrow S_{\l}$ defined by $F_{D,\l}(a)=(D-a)\cap S_{\l}$ and the 
branched $k$-covering $A_{\l}\rightarrow A$ which transforms $F_{D,\l}$ into a single-valued 
holomorphic map from $A_{\l}$ to $S_{\l}$. Clearly, $A_{\l}$ is the ``graph" in $A\times S_{\l}$ 
of the $k$-valued holomorphic map $F_{D,\l}$. 
We thus have a branched covering $\Pi_{\l}:A_{\l}\rightarrow A$. The incidence divisor 
$D_{\l}$ in $A$ is its branch divisor. We fix a non-zero holomorphic 1-form $\om$ on $A$. 
Let $\om_{\l}$ denote its restriction to the curve $S_{\l}$. Let 
$\wt{F}_{D,\l}:A_{\l}\rightarrow S_{\l}$ be the holomorphic map which uniformizes 
the $k$-valued holomorphic map $F_{D,\l}:A\rightarrow S_{\l}$. 
Let $\wt{D}(\om_{\l})_0$ be the pull-back of the canonical divisor $(\om_{\l})_0$ 
of $S_{\l}$ via the holomorphic map $\wt{F}_{D,\l}:A_{\l}\rightarrow A$. 
Let $m:=-(\text{\rm Euler number of $S_{\l}$})=2g(S_{\l})-2=\deg(\om_{\l})$. 
Then $\wt{D}(\om_{\l})_0$ defines a divisor $D(\om_{\l})_0$ on $A$ which is of the form 
$$D(\om_{\l})_0=D_{\l,1}+\cdots+D_{\l,m}$$
where $D_{\l,1},\dots,D_{\l,m}$ are translates of $D$ in $A$. 

\proclaim{Lemma 2.3 (Generalized Hurwitz formula)} 
The incidence divisor $D_{\l} \subset A$, the Gaussian divisor $D_{\l}\subset \PP(TA)$ 
and the given divisor $D\subset A$ fit into the linear equivalence
$$
\pi^*D_{\l}-(\pi^*D_{\l,1}+\cdots+\pi^*D_{\l,m}+D^{[1]}_{\l})=-k\,H^{[1]}\,\,. \tag2
$$
Here, $\pi:\PP(TA) \to A$ is the projection and $H^{[1]}$ is the relative 
hyperplane bundle of the projective vector bundle $\PP(TA) \to A$,  $m=g(S_{\l})$ is 
the genus of the curve $S_{\l}$ and $k=S_{\l}\cdot D$. 
\endproclaim
\noindent
{\it Proof.} Using the local uniformization of the $k$-valued map 
$F_{D,\l}:A\rightarrow S_{\l}$ and the meromorphic vector field $\om_{\l}^{-1}$ 
on $S_{\l}$, we can define a meromorphic section $\Sigma_{\l}$ 
of the relative $\Cal O(k)$ of $\PP(TA)=A\times \PP^{n-1}$ in the following way. 
To $a\in A$ we associate the defining function $\Psi(a)$ 
of the union of $k$ hyperplanes (counted with multiplicities) $\cup_{w\in(D-a)\cap S_{\l}}\PP(T_wD)\subset \PP(T_aA)$ determined by the equation
$$\Psi(a)\,(\underbrace{\om_{\l}^{-1},\dots,\om_{\l}^{-1}}_{k-\text{\rm times}})=1\,\,.$$
This equation determines uniquely the defining function $\Psi(a)$, if all intersections 
of $(D-a)$ and $S_{\l}$ are transversal and $\om_{\l}\not=0$ at these points. 
Since the vector field $\om_{\l}^{-1}$ has poles at points where $\om_{\l}=0$, we have 
extra zeros $D(\om_{\l})_0$ (a divisor on $A$) in $(\Sigma_{\l})_0$ corresponding 
to the zeros of $\om_{\l}$ via the $k$-valued holomorphic map 
$F_{D,\l}:A\rightarrow S_{\l}$. This implies
$$(\Sigma_{\l})_0=D_{\l}^{[1]}+\pi^*D(\om_{\l})_0\,\,.$$
On the other hand, 
if $(D-a)$ and $S_{\l}$ has an intersection point $w$ with multiplicity $\geq 2$, 
the vector $\om_{\l}^{-1}$ at $T_{w}S_{\l}$ belongs to the kernel of the linear form 
defining the tangent space $T_w(D-a)$. 
This implies
$$(\Sigma_{\l})_{\infty}=\pi^*D_{\l}\,\,.$$
Therefore we have the equality of divisors
$$(\Sigma_{\l})_0-(\Sigma_{\l})_{\infty}=D_{\l}^{[1]}+\pi^*D(\om_{\l})_0-\pi^*D_{\l}\,\,.$$
Thus we have the linear equivalence
$$k\,H^{[1]}=D_{\l}^{[1]}+(\pi^*D_{\l,1}+\cdots+\pi^*D_{\l,m})-\pi^*D_{\l}\,\,.$$
This completes the proof of Lemma 2.3.
\qed
\medskip

We note that the right hand side $-k\,H^{[1]}$ of the generalized 
Hurwitz formula does not depend on $D$. 
\medskip

\proclaim{Lemma 2.4} Counting the local intersection number of the projective 
jet lift $f^{[1]}:\C \rightarrow \PP(TA)$ of $f$ with 
the Gaussian divisor $D_{\l}^{[1]} \subset \PP(TA)$ is equivalent to the 
local ramification counting of $f_{\l}:Y_{\l}\rightarrow S_{\l}$. 
Therefore, we have the identification between counting functions
$$N_{f^{[1]},D_{\l}^{[1]}}(r)=N_{f_{\l},\Ram}(r)$$
with a ``slight exception". Here,  the ``slight exception" means the following : 
if $f'(z_0)=0$ holds, the ``order" of the convergence (under the limit 
$w\to z_0$) of $[f'(w)]$ to 
$\bigcup_{w\in (D-f(z_0))\cap S_{\l}}\underbrace{\PP(T_w(D-f(z_0)))}_{\text{\rm translated\,\,to\,\,}f(z_0)}$ may be smaller than expected and thus we may have 
an inequality at the local level in the sense 
that the ramification index of $f_{\l}$ is larger than the intersection number of 
$f^{[1]}$ and $D_{\l}^{[1]}$ at the point under consideration. 
\endproclaim
\noindent
{\it Proof.} This is proved by local consideration and so it suffices to prove under 
the assumption that $f'(z)\not=0$. 
We first consider the case that the intersection of $f^{[1]}:\C\rightarrow \PP(TA)$ 
occurs at $z=z_0\in\C$ such that $f$ touches $D$ at $z=z_0$, i.e., $f(z_0)\in D$. 
Suppose that $D$ is defined by an equation $h(z_1,\dots,z_n)=0$ 
where $z=(z_1,\dots,z_n)$ is the local linear coordinate system of the Abelian variety. 
Let $f:\C \to A$ be a holomorphic curve 
given by $f(z)=(f_1(z),\dots,f_n(z))$. Suppose that near $z=z_0$ we have 
$$h(f_1(z),\dots,f_n(z))=O((z-z_0)^k)\,\,,$$
i.e., the local intersection number of $f$ and $D$ at $z=z_0$ is $k$. 
Differentiating this implies 
$$\frac{\p h}{\p z_1}(f(z))f_1'(z)+\cdots+\frac{\p h}{\p z_n}(f(z))f_n'(z)
=O((z-z_0)^{k-1})\,\,.$$
For each $\l$ such that $f(z_0)+S_{\l}$ is tangent to $D$ at $f(z_0)$ with multiplicity 
$\nu \geq 1$ ($\nu=1$ means the transversal intersection), we define a $\nu$-valued 
holomorphic ``map" $\psi_{\l}$ from a neighborhood of $f(z_0)$ to itself, by assigning 
$w$ (in a neighborhood of $f(z_0)$) the nearest $\nu$ intersection points in 
$(w+S_{\l})\cap D$. Choose a locally defined branch of the $\nu$-valued 
holomorphic map $\psi_{\l}$ and write it using the symbol $\psi_{\l,b}$. 
Then its image is contained in $D$ and therefore $h(\psi_{\l,b}(f(z))=0$ holds for $z$ 
sufficiently close to $z_0$. 
On the other hand, the tangent space at $\psi_{\l,b}(f(z))$ is given by the equation
$$\frac{\p h}{\p z_1}(\psi_{\l,b}(f(z)))\zeta_1
+\cdots+\frac{\p h}{\p z_n}(\psi_{\l,b}(f(z)))\zeta_n=0\,\,.$$
Write $\psi_{\l,b}(f(z))=f(z)+\phi_{\l,b}(z)
=(f_1(z)+\phi_{\l,b,1}(z),\dots,f_n(z)+\phi_{\l,b,n}(z))$. 
Then, from these two equations, we have
$$\split
& \quad \frac{\p h}{\p z_1}(\psi_{\l,b}(f(z)))f_1'(z)+\cdots
+\frac{\p h}{\p z_n}(\psi_{\l,b}(f(z)))f_n'(z)\\
&=\frac{\p h}{\p z_1}(f(z)+\phi_{\l,b}(z))(f_1'(z))+\cdots +\frac{\p h}{\p z_n}(f(z)
+\phi_{\l,b}(f(z))(f_n'(z))\\
&=-\frac{\p h}{\p z_1}(f(z))f_1'(z))+\cdots+\frac{\p h}{\p z_n}(f(z))f_n'(z))
+O((z-z_0)^{\frac{k-1}{\nu}})\\
&=O((z-z_0)^{\frac{k-1}{\nu}})\,\,.
\endsplit
$$
We have $\nu$ such estimates corresponding to the nearest $\nu$ intersection 
points of $S_{\l}$ and $D-f(z_0)$. 
This way, we measure the ``distance" between 
$[f'(z)]\in\PP(T_{f(z)}A)$ and 
$\cup_{w\in (D-f(z))\cap S_{\l}}\underbrace{\PP(T_w(D-f(z)))}_{\text{\rm translated to $f(z)$}}$ 
in terms of the order of divisibility w.r.to $(z-z_0)$. 
Therefore, the local intersection 
number of $f^{[1]}$ and $D_{\l}^{[1]}$ at $z=z_0$ and the local ramification index 
of $f_{\l}$ at $z=z_0$ are the same. We note that this argument does not depend on 
the multiplicity $\nu$ of the intersection of $(D-f(z_0))\cap S_{\l}$ at the neutral point $p$. 

Next, we consider the case that the intersection of 
$f^{[1]}:\C \rightarrow \PP(TA)$ occurs at $z=z_0\in\C$ such that $f(z_0) \not\in D$. 
We get the same local conclusion in this case also. It follows from the definition of 
$D_{\l}^{[1]}$ that we have only to work at such an intersection point of $D-f(z_0)$ 
and $S_{\l}$ with multiplicity $\geq 2$. We want to measure the ``distance" between 
$[f'(z)]\in\PP(T_{f(z)}A)$ and 
$\cup_{w\in (D-f(z))\cap S_{\l}}\underbrace{\PP(T_w(D-f(z)))}_{\text{\rm translated to $f(z)$}}$ 
in terms of the order of divisibility w.r.to $(z-z_0)$. We can do this even if the location 
of $f(z_0)$ is far from $D$. What is important in this consideration is the jet of $f$ at 
$z=z_0$. Indeed, the jet is parallel translated to any place in a portion of $S_{\l}$ 
close to the intersection point with $D-f(z_0)$ by using the group structure 
of the Abelian variety $A$. \qed
\medskip

In the proof of the following Lemma 2.5, we will make a certain observation on the 
counting functions of the ramification of $\pi_{\l}:Y_{\l}\rightarrow \C$ and the intersection 
of $f$ and $D_{\l}$. 
Namely, we will introduce the variation of the Fubini-Study measure of $\L$ (the variation 
parametrized by $t$, positive real numbers in the definition of the counting function) 
which was mentioned after {\it Assumption 1-4} and show that the ``concentration"  
has an effect on the ``averaging" procedure of the difference of two counting functions 
$N_{\pi_{\l},\Ram}(r)-N_{f,D_{\l}}(r)$. 
Namely, we are interested in the ``response" of the average\footnote{\,\,Of course, 
the average over $\L$ depends on the probability measure of $\L$.} 
of two counting functions 
$N_{\pi_{\l},\Ram}(r)$ and $N_{f,D_{\l}}(r)$, under the singular perturbation 
(``concentration") of the measures of $\L$. 
The reason why we are interested in the ``response" is the following. 
The ramification counting function $N_{\pi_{\l},\Ram}(r)$ does not count the order 
of tangency of $f$ to $D$ ``honestly", even if we introduce such $\L$ with large $k$ 
in {\it Assumption 3 and 4}. 
On the other hand, provided we are working under the existence of an upper bound for $r$, 
the counting function $N_{f,D_{\l}}(r)$ counts the tangency of $f$ to $D$ 
``honestly", if we introduce such $\L$ with sufficiently large $k$ (depending on the upper 
bound of $r$) in {\it Assumption 3 and 4}. 

We now introduce the ``variation" of measures on $\L$ in the following way. 
To do so, we take, as an example, the averaging procedure of $N_{f,D_{\l}}(r)$ 
over $\L$ with a variation of measures. Suppose that $f$ and $D$ intersects with 
multiplicity $\geq 2$ at $z=z_0$. If $S_{\l}$ is tangent to $D-f(z_0)$ at the neutral point 
$p$ with sufficiently large intersection multiplicity, this multiple 
intersection point $f(z_0)$ of $f$ and $D$ is counted in the counting function 
$N_{f,D_{\l}}(r)$. This is the reason why we take $N_{f,D_{\l}}(r)$ as an example. 
We recall the definition of the counting function
$$N_{f,D_{\l}}(r)=\int_0^r\frac{n_{f,D_{\l}}(t)-n_{f,D_{\l}}(0)}{t}dt+n_{f,D}(0)\log r\,\,.$$
Fix $r>0$. Let $z_{01},z_{02},\dots,z_{0k}\in\C$ be the places $z$ of $\C$ where there 
exist multiple intersections in $(D-f(z))\cap S_{\l}$ at the neutral point $p$. 
We introduce the following slight modification (if necessary) in counting 
the intersection of $f$ and $D_{\l}$. Namely, if we have several multiple 
intersections of $f$ and $D$ on the circle $\p\C(t)=\{z\in\C\,;\,|z|=t\}$, 
we perturb these points so that these points are separated by the distance 
from the origin, i.e., if $z_{01}$ and $z_{02}$ are such points, i.e., 
$|z_{01}|=|z_{02}|=t$, then we will count these multiple intersections 
(in the definition of $N_{f,D_{\l}}(r)$) as if $|z_{01}|$ and $|z_{02}|$ are very 
close to $t$ and different. Of course we modify only slightly so that the 
difference of the modified counting function and the original counting function 
is negligible in the sense that the difference is at most $\ep\,T_{f,D}(r)$ where 
$\ep>0$ is an arbitrarily given small number. 
We are now ready to introduce the variation of measures on $\L$. 
Let $k$ be a positive number which is sufficiently large compared to the largest 
intersection multiplicity of $f_{\C(r)}:\C(r)\rightarrow A$ and $D$. 
We then introduce the measured space $(\L,d\l_t)$ ($t\leq r$) by the following rule. 
For $t$ such that $0<t< |z_{01}|$, $d\l_t$ is a usual Fubini-Study probability measure. 
Suppose that $D$ is non-singular at the intersection point of $f$ and $D$ 
at $f(z_{01})$. 
For $t$ such that $|z_{0|}|\leq t<|z_{02}|$, we define $d\l_t$ as a regular probability 
measure ``concentrated" at $\L_{j_1}$ corresponding to the hyperplane in $T_pA$ which is 
the parallel translate to the neutral point $p$ of the tangent plane of $D$ at 
the multiple intersection point $f(z_{01})$ 
with $f$ (here, we assume that $D$ is non-singular at $f(z_{01})$), 
where $j$ is the intersection multiplicity of $f$ and $D$ at $z=z_{01}$. 
If $f(z_0)$ belongs to the singular locus of $D$, we need a slight modification. 
However, the modification is easy, because $D$ is assumed to have at worst normal 
crossings. 
Next, suppose that $D$ is non-singular at the intersection point of $f$ and $D$ 
at $f(z_{02})$. For $t$ such that $|z_{02}|\leq t<|z_{03}|$, we define $d\l_t$ as 
a regular probability measure ``concentrated" at $\L_{j_2}$ corresponding to 
the hyperplane in $T_pA$ which is the parallel translate to the neutral point $p$ 
of the tangent plane of $D$ at the multiple intersection point $f(z_{02})$ 
with $f$ (here, we assume that $D$ is non-singular at $f(z_{02})$), 
where $j_2$ is the intersection multiplicity of $f$ and $D$ at $z=z_{02}$. Again 
the necessaey modification when $D$ is singular at $f(z_{02})$ is easy. 
Repeating this procedure, we get a variation of probability measures 
$d\l_t$ ($0<t<r$) and get the measured space $(\L,d\l_t)$ parametrized by $t$ 
($0<t<r$). We then define the ``average of the counting function $N_{f,D_{\l}}(r)$ 
over the variation of  probability spaces $(\L,d\l_t)$", denoted by 
$\displaystyle N^{(\L,d\l_t)}_{f,D_{\l}}(r)$, as follows:
$$
N^{(\L,d\l_t)}_{f,D_{\l}}(r)
:=\int_0^r\frac{dt}{t}\int_{(\L,d\l_t)}(n_{f,D_{\l}}(t)-n_{f,D_{\l}}(0))d\l_t+n_{f,D}(0)\log r\,\,.
$$
For other Nevanlinna theoretic functions, we define the corresponding ``average" 
over the variation of  probability spaces $(\L,\l_t)$ in the following way:
$$
 \split
& T^{(\L,d\l_t)}_{f,D_{\l}}(r)
:=\int_0^r\frac{dt}{t}\int_{(\L,d\l_t)}d\l_t\int_{\C(t)}f^*c_1(\Cal O_A(D_{\l}))\,\,,\\
& m^{(\L,d\l_t)}_{f,D_{\l}}(r)
:=\int_0^{2\pi}\frac{d\theta}{2\pi}\int_{(\L,d\l_r)}\log\frac1{\dist_{\Euc}(f(re^{i\theta}),D_{\l})}\,\,,
\endsplit
$$
The response under this variation of probability spaces is summarized 
in the following:

\proclaim{Lemma 2.5} For any variation of probability spaces $(\L,d\l_t)$ ($0<t<r$) 
as above, we have
$$
m^{(\L,d\l_t)}_{f_{\l},p}(r)\,=\,m_{f,D}(r)+O(1)\,\,.\tag3
$$
On the other hand, let $\ep$ be any positive number. Then there exists 
a variation of probability measures $(\L,d\l_t)$ ($0<t<r$) with sufficiently 
strong ``concentration" such that the ``averaging" over $(\L,d\l_t)$ 
satisfies the following estimates : 
$$
\align
& N^{(\L,d\l_t)}_{f_{\l},\Ram}(r)-N^{(\L,d\l_t)}_{f^{[1]},D_{\l}^{[1]}}(r)\,=0\,\,,\tag4\\
& N^{(\L,d\l_t)}_{\pi_{\l},\Ram}(r)-N^{(\L,d\l_t)}_{f,D_{\l}}(r) \,\leq\, -N^1_{f,D}(r)
+\ep\,T_{f,D}(r)\,\,.\tag5
\endalign
$$
Here, $p$ in (3) is the neutral point (see the beginning of \S2), the 
$O(1)$ term in (3) is independent of the variation $(\L,d\l_t)$ 
and $N^1_{f,D}(r)$ in (4) 
is the residual counting function counting function $N_{f,D}(r)-N_{1,f,D}(r)$, 
where $N_{1,f,D}(r)$ is the counting function which counts all intersections 
with multiplicity $1$. Moreover, the equality (4) holds ``with slight exceptions" 
(i.e., in the local level the equality $=$ may be inequality $\geq$ at $z$ 
such that $f'(z)=0$) in the same sense as in Lemma 2.4. \endproclaim
\noindent
{\it Proof.} Since $p\in S_{\l}$ for all $\l\in\L$ ($p$ is the neutral point of $A$) 
and $f$ approximates $D$ if and only if $f_{\l}$ approximates $p$, 
we have the averaging formula (3). 
Indeed, it is easy to prove (3) by combining the argument in Lemma 2.4 and 
[K3, Lemma 2.2]. 
Indeed, [K3, Lemma 2.2] implies that the proof of (3) is based on the estimate 
of the functions on the projective space which has logarithmic poles along the 
hyperplane at infinity and Lemma 2.4 implies that the estimate is stable for 
large $k$ with bound (recall that, in the proof of Lemma 2.4, we have the multiplicity 
$\nu$ in the denominator of the exponent of $(z-z_0)$ and we have $\nu$ such 
estimates). The estimate ``with slight exception'' is a consequence of Lemma 2.4. 
To prove (5), we choose the family of curves $\{S_{\l}\}_{\l\in\L}$ with $k$ in 
{\it Assumption 3 and 4} sufficiently large compared to the maximum multiplicity 
of $f$ and $D$ in $\C(r)$. 
We can then choose a variation of probability spaces $(\L,\l_t)$ in such a way that it 
picks up all multiple intersections of $f$ and $D$ in $\C(r)$ just like in the above 
argument. 
If the intersection multiplicity of $S_{\l}$ and $D-f(z_0)$ at $p$ is sufficiently large 
compared to the intersection multiplicity of $f$ and $D$ at $z=z_0$, then 
the intersection of $f$ and $D_{\l}$ is equivalent to that of $f$ and $D$ at $z=z_0$. 
In this situation, the counting function $N_{f,D_{\l}}(r)$ counts all multiple intersections 
of multiplicity $\nu$ (of $f$ and $D$ in $\C(r)$) as those of multiplicity $\nu-1$. 
On the other hand, the counting function $N_{\pi_{\l},\Ram}(r)$ does not count the 
multiplicity of $f$ and $D$. Indeed, the $k$-covering map $\pi_{\l}:Y_{\l}\rightarrow\C$ 
is obtained by the pull back (via $f:\C\rightarrow A$) of the branched $k$-covering 
$A_{\l}\rightarrow A$ uniformizing the $k$-valued holomorphic map 
$A\rightarrow S_{\l}$ defined by $A\ni a \mapsto (D-a)\cap S_{\l}$ (cf. the proof of 
Lemma 2.3). Therefore, given small $\ep>0$, if the ``concentration" 
to the sequence of $\L_j$'s which appear corresponding to the sequence 
of multiple intersection points (of $f$ and $D$ in $\C(r)$) is sufficiently strong, 
then we have the estimate of the form of (5).  \qed
\medskip

The following corollary is a consequence of (3) in Lemma 2.5 (computation in the 
proof of Lemma 2.4), or more generally, the invariance 
of the proximity functions under ``measure concentration". 

\proclaim{Corollary 2.6} The First Main Theorem 
$m_{f,D_{\l}}(r)+N_{f,D_{\l}}(r)=T_{f,D_{\l}}(r)+O(1)$ remains to be true after taking the 
``average" over the variation of probability spaces 
$(\L,d\l_t)$, i.e., we have
$$m^{(L,d\l_t)}_{f,D_{\l}}(r)+N^{(\L,d\l_t)}_{f,D_{\l}}(r)=T^{(\L,d\l_t)}_{f,D_{\l}}(r)+O(1)\,\,.$$
\endproclaim

To sum up, we have introduced the Gaussian divisor $D_{\l}^{[1]}$ and the 
incidence divisor $D_{\l}$ to interpret the ramification counting functions 
$N_{f_{\l},\Ram}(r)$ and $N_{\pi_{\l},\Ram}(r)$. 
To get information which is independent of $\l$, we have to take the ``average" 
over $\L$ against some probability measure. 
To get such information, we have introduced the ``variation of probability 
spaces" $(\L,d\l_t)$ and look at the difference of the response under this 
``variation of probability measures" of $\L$. Namely,  the ``concentration" of 
probability measures at subvarieties $\L_j$'s of $\L$ corresponding to high 
degree of tangents of $S_{\l}$ and $D$. 
As a result, the difference $N_{f_{\l},\Ram}(r)-N_{f^{[1]},D_{\l}^{[1]}}(r)$ 
is essentially stable under the ``variation of probability space" $(\L,d\l_t)$. 
On the contrary, the quantity $-N^1_{f,D}(r)$ stemming from the 
multiple intersection (of $f$ and $D$ in $\C(r)$), which is only a small portion 
of the difference $N_{\pi_{\l},\Ram}(r)-N_{f,D_{\l}}(r)$ w.r.to the average against the 
usual Fubini-Study measure, becomes ``dominating" under the average against 
the variation of the probability space $(\L,d\l_t)$ with sufficiently strong ``concentration". 
\medskip

There is a natural way\footnote{\,\, We use the stratification 
of $\L$ into $\L_j$'s w.r.to a fixed hyperplane of $T_pA$ in {\it Assumption 4}.} 
of producing the variation of probability spaces 
$(\L,d\l_t)$ which is useful in the proof of Lemma 2.5. 
We introduce the notion of the ``combinatorial blow up" of $\L$. 
We ``expand" the parameter space $\L$ in a combinatorial way 
without changing the family of curves itself. 
Here, expanding the parameter space $\L$ in a combinatorial way {\it e.t.c.} means 
the following. 
Suppose that $\l$ is a positive number valued 
coordinate function of $\L$ (for instance, some angle parameter naturally defined on $\L$ 
by using the Hopf fibration from odd dimensional sphere to a complex projective space). 
If we consider the totality of all decompositions of $\l$ into $N$ 
positive numbers $\l=\l_1+\cdots+\L_N$, we can increase the degree of 
freedom in the coordinate\footnote{\,\,This is something 
similar to the blow up, because we replace a point $\l_0$ on a $\l$-line by 
the $(N-1)$-simplex $\{(\l_1,\dots,\l_N)\,|\,\l_1+\cdots+\L_N=\l_0,\,\l_1,\dots,\l_N>0\}$.}. 
 We call this the ``combinatorial blow up". We then introduce the natural probability 
measure on the combinatorial blow up. This has the effect that the neighborhood of 
the fixed radius of a point becomes relatively smaller, if $N$ becomes larger. 
Applying  the combinatorial blow up to $\L$, we can increase its dimension 
without increasing the family itself (i.e., this increases only the parameter space 
without increasing the family). 
The tuple $(\l_1,\dots,\l_N)$ represents the same $S_{\l}$ for $\forall (\l_N^1,\dots,\l_N^N)$ 
if $\l=\l_N=\l_N^1+\cdots+\l_N^N$.
\medskip

The following Lemma 2.7 is a direct consequence of Lemma 2.5 and explains the 
reason why one can reduce the truncation level to $1$ in the Second Main Theorem 
for holomorphic curves into Abelian varieties. 
Let $N^{(\L,d\l_t)}_{f,D_{\l}}(r)$, {\it etc}, be the same as in Lemma 2.5. 

\proclaim{Lemma 2.7} Let $A$ be an Abelian variety, $D$ any ample reduced divisor, 
$\ep>0$ any positive number and $f:\C \to A$ any holomorphic curve such that 
$f(\C) \not\subset \Supp(D)$. 
Then there exists a variation of probability measures $(\L,d\l_t)$ ($0<t<r$) 
with sufficiently strong ``concentration" such that the ``averaging" over $(\L,d\l_t)$ 
satisfies the following estimates : 
$$
\split
& N^{(\L,d\l_t)}_{f_{\l},\Ram}(r)
-N^{(\L,d\l_t)}_{f^{[1]},D_{\l}^{[1]}}(r)+N^{(\L,d\l_t)}_{f,D_{\l_N}}(r)-N^{(\L,d\l_t)}_{\pi_{\l},\Ram}(r)\\
& \quad \geq N^1_{f,D}(r)-\ep\,T_{f,D}(r)\exc\,\,.
\endsplit
$$
\endproclaim

\comment
\noindent
{\it Proof.} The inequality (5) means the following. 
Let $N'_{f^{[1]},D_{\l}^{[1]}}(r)$ be the counting function which counts the intersection 
of $f^{[1]}:\C\to A^{[1]}$ and $D^{[1]}_{\l}\subset A^{[1]}$ only at $z$ such that 
$\deg_z(f^*D)\geq 2$. Then we have the averaging formula:
$$
\int_{\l_N\in\L_N}N'_{f^{[1]},D_{\l_N}^{[1]}}(r)\,d\l_N \leq \,\ep\,T_{f,E}(r) \exc\,\,. \tag6
$$

Next, we compare the ramification of $\pi_{\l_N}:Y_{\l_N} 
\rightarrow \C$ to the intersection multiplicity of the holomorphic 
curve $f:\C \rightarrow A$ and the incidence divisor $D_{\l_N}$. For this purpose 
extending $D$ is not necessary. 
The intersection multiplicity of the holomorphic curve $f:\C\rightarrow A$ 
with the incidence divisor $D_{\l_N}$ is $\geq 2$, if the ramification index 
of $\pi_{\l_N}$ is $\geq 1$ (we note that $Y_{\l_N}$ is the normalization of 
the covering surface over $\C$ defined by collecting all intersections 
of $(D-f(z))\cap S_{\l_N}$ over $z\in \C$). More precisely, the difference is 
always nonnegative. Indeed, if $f(z) \in D_{\text{\rm reg}}$, 
then the intersection index of $f$ and $D_{\l_N}$ at $f(z)$ and the ramification 
index of $\pi_{\l_N}:Y_{\l_N} \to \C$ is nonnegative for almost all $\l_N$. 
If $f(z)\in \Sing(D)$, then the intersection $(D-f(z))\cap S_{\l_N}$ ($\forall \l_N$) 
always contains the neutral point as a multiple point. 
Therefore, the intersection multiplicity of $f$ and $D$ has an 
effect on the ramification index of $f_{\l_N}$. 
We show that the difference of the intersection multiplicity of $f$ and $D_{\l_N}$ and 
the ramification index of $\pi_{\l_N}:Y_{\l_N} \rightarrow \C$ is nonnegative. 
Indeed, there exists a Zariski open subset $(\L_N)_{z_0}$ of $\l_N$ such that 
for any $\l_N \in (\L_N)_{z_0}$ and for any $z\not= z_0$ close to $z_0$, the curve 
$S_{\l_N}+f(z)$ and $D$ intersect at distinct $\deg_{z_0}(f^*D)$ points near $f(z_0)$. 
This implies that $N_{f,D_{\l_N}}(r) \geq N_{\pi_{\l_N},\Ram}(r)$ holds in the local level. 
Therefore we have
$$
\int_{\l_N\in\l_N}\{N_{f,D_{\l_N}}(r)-N_{\pi_{\l_N},\Ram}(r)\}d\l_N \geq 0\,\,.\tag7
$$
Combining (6) and (7), we obtain the desired inequality of Lemma 2.5. 
Indeed, we can argue as follows. The definition of the ramification counting 
function implies
$$
N_{f_{\l_N},\Ram}(r)=N_{f,D}(r)-N_{1,f,D}(r)\,\,(:=N^1_{f,D}(r))\,\,,
$$
where $N_{1,f,D}(r)$ is the usual truncated counting function with truncation level 
$1$, i.e., the counting function ``without counting multiplicity".  
Therefore, we have
$$
\split
& \quad \int_{\l_N\in\l_N}\{N_{f_{\l_N},\Ram}(r)-N_{f^{[1]},D_{\l_N}^{[1]}}(r)+N_{f,D_{\l_N}}(r)
-N_{\pi_{\l_N},\Ram}(r)\}\,d\l_N\\
& \geq \int_{\l_N\in\l_N}\{N'_{f_{\l_N},\Ram}(r)-N'_{f^{[1]},D_{\l_N}^{[1]}}\}\,d\l_N\\
& \geq (N_{f,D}(r)-N_{1,f,D}(r))-\ep\,T_{f,E}(r)\exc\\
& =N^1_{f,D}(r)-\ep\,T_{f,E}(r)\exc \,\,.
\endsplit
$$
\endcomment

We are ready to prove our main result (1). 

Let $f:\C\rightarrow A$ be a holomorphic curve into an Abelian variety $A$ and $D$ 
an ample reduced divisor in $A$ (for simplicity, we assume that $D$ has at worst normal 
crossings). 

From here on, we consider the family of algebraically equivalent algebraic curves 
$\{S_{\l}\}_{\l\in\L}$ satisfying the {\it Assumptions 1-4} stated at the beginning of this 
\S2. In particular, the parameter space $\L$ is equipped with the variation 
of the probability measures $d\l_t$ so that Lemma 2.5 and 2.7 hold. 
In particular, we consider the variation of the probability spaces $(\L,d\l_t)$ 
where the measure $d\l_t$ ``concentrates" at $\L_j$'s sufficiently strongly, 
so that we have the estimates in (3), (4), (5) in Lemma 2.5 and the estimate in 
Lemma 2.7. 
\medskip

We start by applying the Second Main Theorem to the family of 
holomorphic maps $f_{\l}:Y_{\l} \rightarrow \C$ and ``average" the result 
over the variation of probability spaces $(\L,\l_t)$.  

The Second Main Theorem for holomorphic maps from a finite analytic covering 
space of $\C$ to a compact Riemann surface is available in, for instance, [N]. 
Let $p$ denote the neutral point of the abelian variety $A$. 
Then for all $\l\in\L$ the curve $S_{\l}$ contains $p$. 
Now the Second Main Theorem for $f_{\l}:Y_{\l} \rightarrow \C$ is stated 
as follows:
$$
m_{f_{\l},p}(r)+T_{f,K_{S_{\l}}}(r)+N_{f_{\l},\Ram}(r)
-N_{\pi_{\l},\Ram}(r) \leq O(\log T_{f,D}(r)+\log r)\exc\,\,.
\tag8
$$
The right hand side is bounded above by the quantity of the form $k\,\log T_{f,D}(r)$, where 
$k=S_{\l}\cdot D=[Y_{\l}:\C]$ is the covering degree of $\pi_{\l}:Y_{\l}\rightarrow \C$. 

Since $p\in S_{\l}$ for all $\l\in\L$ and $f$ approximates $D$ if and only 
if $f_{\l}$ approximates $p$, (3) in Lemma 2.5 implies the averaging formula
$$
m_{f,D}(r)=m^{(\L,d\l_t)}_{f_{\l},p}(r)d\l+O(1)\,\,,
$$
where the $O(1)$-term is independent of $f$. 

Next we would like to apply the above Second Main Theorem to $f_{\l}$ and 
average the result over the variation of probability spaces $(\L,\l_t)$ in the sense of 
Lemma 2.5. 
Here, the existence of exceptional intervals in the Second Main Theorem 
may cause difficulty. 
However, the size of the exceptional intervals implicit in the symbol $\exc_{\ep}$ 
is controlled by the Borel Lemma, whose ``family version" is formulated and 
proved in [K2, (2.6)]. 
This implies that we can average the Second Main Theorems for $f_{\l}$'s 
without worrying about the individual exceptional intervals. Moreover, 
the Second Main Theorem is a consequence of Nevanlinna's lemma on 
logarithmic derivative (see, for instance [NoO] and [L2]). 
The constant term depends on $\l$ 
and might cause problem in the course of averaging. This constant term in turn 
stems from the constant term in the First Main Theorem for $f_{\l}$ and is 
of the form $\log|c_{\l}|$, 
where $c_{\l}$ is determined by the behavior of $f_{\l}$ at $z=0$. 
As $f$ is fixed and $f_{\l}$ is defined algebraically using 
the intersection $(D-f(z))\cap S_{\l}$ ($\l\in\L$), $c_{\l}$ depends 
algebraically on $\l$ and therefore the integral 
of $\log|c_{\l}|$ over $\L$ is finite whose value 
depends only on the behavior of $f$ at $z=0$ and $\{S_{\l}\}_{\l\in\L}$. 

We may now apply the Second Main Theorem (8) to each $f_{\l}$ and average 
the result over $(\L,\l_t)$ in the sense of Lemma 2.5. Using the concavity of the 
logarithm, we have
$$
\split
& \quad m_{f,D}(r) = m^{(\L,d\l_t)}_{f_{\l},p}(r)d\l+O(1)\\
& \leq -T^{(\L,d\l_t)}_{f_{\l},K_{S_{\l}}}(r)-(N^{(\L,d\l_t)}_{f_{\l},\Ram}(r)
-N^{(L,d\l_t)}_{\pi_{\l},\Ram}(r))
+O(\log T^{(\L,d\l_t)}_{f_{\l},p}(r)+\log r)\exc\,\,.
\endsplit
$$
It follows from the definition that the averaging formula
$$
T_{f,D}(r)=T^{(\L,d\l_t)}_{f_{\l},p}(r)d\l+O(1)
$$
holds, where $O(1)$-term does not depend on $f$. Let $\deg K_{S_{\l}}=m$ and 
$D_{\l,1},\dots,D_{\l,m}$ be as in Lemma 2.3. Then, we have
$$
T^{(\L,d\l_t)}_{f_{\l},K_{S_{\l}}}(r)d\l=T_{f,D_{\l,1}+\cdots+D_{\l,m}}(r)+O(1)\,\,.
$$
Applying Lemma 2.7, we get
$$
\split
& m_{f,D}(r)\leq \{-T^{(\L,d\l_t)}_{f,D_{\l,1}+\cdots+D_{\l,m}}(r)-N^{(\L,d\l_t)}_{f^{[1]},D_{\l}^{[1]}}(r)
+N^{(\L,d\l_t)}_{f,D_{\l}}(r)\}\}\\
& \quad -\{(N^{(\L,d\l_t)}_{f_{\l},\Ram}(r)
-N^{(\L,d\l_t)}_{f^{[1]},D_{\l}^{[1]}}(r))-(N^{(\L,d\l_t)}_{\pi_{\l},\Ram}(r)-N^{(\L,\l_t)}_{f,D_{\l}}(r))\}\\
& \quad +O(\log T_{f,D}(r)+\log r)+O(1)\exc\,\,\\
& \leq \{-T^{(\L,d\l_t)}_{f,D_{\l,1}+\cdots+D_{\l,m}}(r)-N_{f^{[1]},D_{\l}^{[1]}}(r)
+N_{f,D_{\l}}(r)\}\\
& \quad -N^1_{f,D}(r)
+\ep\,T_{f,D}(r)+O(\log T_{f,D}(r)+\log r)+O(1)\exc\,\,.
\endsplit
$$
We now apply the First Main Theorem (Corollary 2.6) to replace 
the counting function $N_{f^{[1]},D_{\l}^{[1]}}(r)$ (resp. $N_{f,D_{\l}}(r)$) 
by the difference of the height function and the proximity 
function $T_{f^{[1]},D_{\l}^{[1]}}(r)-m_{f^{[1]},D_{\l}^{[1]}}(r)$ 
(resp. $T_{f^{[1]},\pi^*D_{\l}}(r)-m_{f^{[1]},\pi^*D_{\l}}(r)$):
$$
\split
m_{f,D}(r)+N^1_{f,D}(r) =  & \{-T^{(\L,d\l_t)}_{f^{[1]},\pi^*D_{\l,1}+\cdots+\pi^*D_{\l,m}}(r)
-T^{(\L,d\l_t)}_{f^{[1]},D_{\l}^{[1]}}(r)+T^{(\L,d\l_t)}_{f^{[1]},\pi^*D_{\l}}(r)\}\\
& +m^{(\L,\l_t)}_{f^{[1]},D_{\l}^{[1]}}(r)
-m^{(\L,d\l_t)}_{f^{[1]},\pi^*D_{\l}}(r)\\
& +\ep\,T_{f,D}(r)+O(\log T_{f,D}(r)+\log r)+O(1)\exc\,\,.
\endsplit
$$
As $\{D_{\l}\}_{\l\in\L}$ is free from base locus, the computation in the proof 
of Lemma 2.5 implies
$$
\int_{\l\in\L}m_{f^{[1]},\pi^*D_{\l}}(r)d\l=O(1)\,\,.
$$
Hence
$$
\split
m_{f,D}(r)+N^1_{f,D}(r) \,\, \leq\,\, 
& T^{(\L,d\l_t)}_{f^{[1]},-m\,\pi^*D-D_{\l}^{[1]}+\pi^*D_{\l}}(r)
+m^{(\L,d\l_t)}_{f^{[1]},D_{\l}^{[1]}}\\
& +\ep\,T_{f,D}(r)+O(\log T_{f,D}(r)+\log r)+O(1)\exc\,\,.
\endsplit
$$
Applying the Generalized Hurwitz formula (2) in Lemma 2.3, we have
$$
\split
m_{f,D}(r)+N^1_{f,D}(r) \,\,\leq\,\, & -k\,T_{f^{[1]},H^{[1]}}(r)
+m^{(\L,d\l_t)}_{f^{[1]},D^{[1]}_{\l}}(r)\\
& +\ep\,T_{f,D}(r)+O(\log T_{f,D}(r)+\log r)+O(1)\exc\,\,.
\endsplit
$$
It follows from the definition of the Gaussian divisor that
$$
\split
\text{\rm Bs}\,\{D_{\l}^{[1]}\}_{\l\in\L} & 
=\text{\rm Zariski\,\,closure\,\,of\,\,}\PP(TD_{\text{\rm reg}})\\
& =:D^{[1]}\quad (\subset \PP(TA))\,\,.
\endsplit
$$
This implies that, averaging $m_{f^{[1]},D_{\lambda}^{[1]}}(r)$ 
over $\L$, we have a main contribution of the base locus $D^{[1]}$ 
as well as that coming from those curves $S_{\lambda}$ 
which are tangent to $D-w$ ($w\in D$) at the origin. The latter 
is not larger than $\ep\,T_{f,D}(r)$ for any positive $\ep$. 
Namely we have the following Lemma. 

\proclaim{Lemma 2.8} We have
$$
m^{(\L,d\l_t)}_{f^{[1]},D^{[1]}_{\l}}(r)d\l\,\, 
\leq \,\, m_{f^{[1]},D^{[1]}}(r) + \ep\,T_{f,D}(r)\,\,.
$$
\endproclaim
\noindent
{\it Proof.} We write 
the average $m^{(\L,d\l_t)}_{f^{[1]},D_{\lambda}^{[1]}}(r)$ as 
$$
\int_0^{2\pi}\biggl(\int_{(\Lambda,\l_r)}
\log^+\frac1{\dist_{\Euc}(f^{[1]}(re^{i\theta}),D_{\lambda}^{[1]})}
d\lambda_r\biggr)\frac{d\theta}{2\pi}\,\,.
$$
The main part of the integral over $\Lambda$ consists of the contribution 
of the base locus $D^{[1]}$, i.e., 
$$
\log^+\frac1{\dist_{\Euc}(f^{[1]}(re^{i\theta}),D^{[1]})}
$$
and possibly the contribution from those curves $S_{\lambda}$ 
which are tangent to $D-w$ ($w\in D$) at the neutral point. 
To see the reason why we have to take the latter contribution into 
account, we take an arbitrary $w\in D$. 
Let $\l\in\L$ be such that $S_{\l}$ is tangent to $D-w$ at the 
neutral point. Then the Gaussian divisor $D_{\l}^{[1]}$ contains $D^{[1]}$ 
with multiplicity $\geq 2$. So, if $f$ approximates $w\in D$, we have to 
take special care of those $S_{\l}$'s tangent to $D-w$ at the neutral point. 
The latter contribution is not larger than the product of the 
volume of the $\ep$-tubular neighborhood of the subvariety $\L_w$ 
of $\L$ parameterizing those curves and the proximity 
function of $f^{[1]}$ to the Gaussian divisor $D_{\lambda}^{[1]} \subset 
\PP(TA)$ ($\l\in\L_w$). 
So, if the measure under consideration is ``concentrated", the main contribution 
comes from those $\l$ such that the corresponding $S_{\l}$ intersects $D-w$ 
($w\in D$) at the neutral point with high multiplicity. 
However, it follows from the  computation in Lemma 2.4 that the contribution 
to the proximity function $m_{f^{[1]},D_{\l}^{[1]}}(r)$ is independent of the 
multiplicity of $S_{\l}$ and $D-w$. Therefore, the main part of the left hand 
side of Lemma 2.6 is $m_{f^{[1]},D^{[1]}}(r)$. 
The error of this approximation comes from the portion of $D_{\l}^{[1]}$ 
close to the base locus $D^{[1]}$. 
On the other hand, the approximation of $f^{[1]}$ to $D_{\lambda}^{[1]}$ 
is dominated by the height function $T_{f^{[1]},D_{\lambda}^{[1]}}(r)$ and 
this height function is bounded from above and below by the 
height function $T_{f,D}(r)$ in the sense that there exist positive 
constants $C_1$ and $C_2$ such that
$$
C_1\,T_{f,D}(r) \leq T_{f^{[1]},D_{\l}^{[1]}}(r) 
\leq C_2\,T_{f,D}(r)
$$
holds. Therefore the the latter contribution is not larger 
than $\ep\,T_{f,D}(r)$ for any positive number $\ep$ in the asymptotic sense 
when $r \to \infty$. \qed
\medskip

Summing up the above argument, we have

\proclaim{Proposition 2.9} Let $A$ be an Abelian variety, $D$ a reduced 
divisor and $\{S_{\l}\}_{\l\in\L}$ an algebraic family of algebraically 
equivalent curves in $A$ all passing through the neutral point $p$ such 
that general $S_{\l}$ is a non-singular curve and the natural rational 
map $\L \mer \PP(T_pA)$ is surjective. 
Let $H^{[1]}$ (resp. $D_{\l}^{[1]}$) denote the relative hyperplane 
bundle of $\PP(TA)$ (resp. the Gaussian divisor with respect 
to $\{S_{\l}\}_{\l\in\L}$). Let $\ep$ be any positive number. 
Then we have
$$
\split
m_{f,D}(r)+N^1_{f,D}(r)\,\,
\leq \,\, & -k\,T_{f^{[1]},H^{[1]}}(r)+m_{f^{[1]},D^{[1]}}(r)\\
& +\ep\,T_{f,D}(r)+O(\log T_{f,D}(r)+\log r)+O(1)\exc_{\ep}\,\,.
\endsplit
$$
\endproclaim

So far we worked over $A$ and $A^{[1]}=\PP(TA)$. However, we are able to 
perform all above arguments over the cone $A^{(1)}:=\PP(TA\oplus \Cal O_A)$. 
Here we list all definitions which should be modified if we replace $A^{[1]}$ 
with $A^{(1)}$. 
\NI
(i) The Gaussian divisor $D_{\l}^{[1]}$ should be replaced with its 
cone $D_{\l}^{(1)}$ (i.e., $A^{[1]}$ is the hyperplane at infinity in $A^{(1)}$ 
and $D^{(1)}$ is defined to be the cone over $D^{[1]}_{\l}$ in $A^{(1)}$). 
\NI
(ii) The projective jet lift $f^{[1]}$ should be replaced with its 
cone $f^{(1)}(z):=[f'(z)\oplus 1_{f(z)}]$, where $1_a$ ($a\in A$) is the value of 
the generator of $H^0(A,\Cal O_A)$ at $a$. 
\NI
(iii) The relative hyperplane bundle $H^{[1]}$ of the projective vector 
bundle $A^{[1]}\to A$ should be replaced with the relative hyperplane bundle $H^{(1)}$ 
of the cone $A^{(1)} \to A$. 
\NI
(iv) $D^{[1]}$ should be replaced with its cone $D^{(1)}$.

All the other definitions ($D_{\l}$, $Y_{\l}$, $\pi_{\l}:Y_{\l}\to\C$ 
and $f_{\l}:Y_{\l}\to S_{\l}$) are the same. 

The Generalized Hurwitz formula (2) remains true, if we replace $H^{[1]}$ 
(resp. $\pi:A^{[1]}\to A$ and $D^{[1]}_{\l}$) with $H^{(1)}$ 
(resp. their cone $\pi:A^{(1)} \to A$ and $D^{(1)}_{\l}$). 

The crucial Lemmas 2.5 and 2.7 continue to hold, if we replace $f^{[1]}$ 
(resp. $D^{[1]}_{\l}$) with its ``cone'' $f^{(1)}$ (resp. $D^{(1)}_{\l}$). 

As a result we have

\proclaim{Proposition 2.10} Let $A$ be an Abelian variety, $D$ a reduced 
divisor and $\{S_{\l}\}_{\l\in\L}$ an algebraic family of algebraically 
equivalent curves in $A$ all passing through the neutral point $p$ such 
that general $S_{\l}$ is a non-singular curve and the natural rational 
map $\L \mer \PP(T_pA)$ is surjective. 
Let $H^{(1)}$ (resp. $D_{\l}^{(1)}$) denote the relative hyperplane 
bundle of $\PP(TA\oplus \Cal O_A)$ (resp. the cone of the Gaussian divisor 
with respect to $\{S_{\l}\}_{\l\in\L}$). Let $\ep$ be any positive number. 
Then we have
$$
\split
m_{f,D}(r)+N^1_{f,D}(r)\,\,
\leq \,\, & -k\,T_{f^{(1)},H^{(1)}}(r)+m_{f^{(1)},D^{(1)}}(r)\\
& +\ep\,T_{f,D}(r)
+O(\log T_{f,D}(r)+\log r)+O(1)\exc_{\ep}\,\,.
\endsplit
$$
\endproclaim

Therefore, the proof of (1) is reduced to the proof of the estimate 
of type
$$m_{f^{(1)},D^{(1)}}(r)\leq \ep\, T_{f,D}(r)\,\,.$$
In \S3, we solve this problem by using a geometric version of 
Nevanlinna's Lemma on Logarithmic Derivative ([Y1], [K3]). 

\beginsection 3. Lemma on logarithmic derivative.

It follows from Proposition 2.10 that the proof of (1) is reduced to the 
estimates of $T_{f^{(1)},H^{(1)}}(r)$ and $m_{f^{(1)},D^{(1)}}(r)$. 
For this purpose we consider the sequence of successive projective completion 
of tangent bundles
$$
\split
& A^{(0)}=A,\,\,A^{(1)}=\PP(TA\oplus \Cal O_A),\,\,
A^{(2)}=\PP(TA^{(1)}\oplus \Cal O_{A^{(1)}}),\dots,\\
& A^{(i)}=\PP(TA^{(i-1)}\oplus \Cal O_{A^{(i-1)}}),\dots\,\,.
\endsplit
$$
A geometric version of the 
Lemma on logarithmic derivative ([Y1], [V2] and [K1]) 
yields the estimates for $T_{f^{(1)},H^{(1)}}(r)$ 
and $m_{f^{(1)},D^{(1)}}(r)$. 
Here we recall this. Let $X$ be a smooth projective variety. 
Write $X^{(i)}$ for the $i$-th successive projective completion of 
tangent bundles. Let $Z$ a 
subscheme given by $Z=V(f_1,\dots,f_k)$. 
Then $Z^{(i)}$ is a subvariety of $X^{(i)}$ defined by the Zariski closure 
of $Z_{\text{\rm reg}}^{(i)}$ in $A^{(i)}$. By the symbol $\infty$, we denote 
the divisor at infinity in the projective completion of a vector bundle. 
We are now ready to state a geometric version of Lemma on logarithmic 
derivative. 

\proclaim{Lemma 3.1} Let $X$ be a smooth projective variety, $Z$ any subscheme 
and $E \to X$ any ample line bundle. Let $f:\C \to X$ be an arbitrary 
holomorphic curve such that $f(\C) \not\subset \text{\rm Supp}(Z)$. 
Then we have
$$
\left\{\aligned
& m_{f^{(i)},Z^{(i)}}(r) \leq m_{f,Z}(r)+O(\log^+(r\,T_{f,E}(r)))\exc\,\,,\\
& m_{f^{(i)},\infty}(r) \leq O(\log^+(r\,T_{f,E}(r)))\exc\,\,.
\endaligned\right.
$$
\endproclaim

The advantage of this formulation is that the original form 
$$
m_{f'/f,\infty}(r) 
\leq O(\log^+(rT_f(r)))\exc
$$
of Nevanlinna's lemma on logarithmic derivative for $f:\C \to \PP^{1}(\C)$ 
splits into two inequalities each of which has its clear geometric meaning. 
For the proof, we refer to [Y1], [V2] and [K1,3]. 

We return to our situation. We apply Lemma 3.1 to the holomorphic 
curve $f:\C \to A$. We put $X=A$, $Z=D$ in Lemma 3.1. We fix an ample 
line bundle $E \to A$. 

The second inequality of Lemma 3.1 implies
$$
\split
T_{f^{(1)},H^{(1)}}(r) & = \int_0^{2\pi}\log(1+||f'(re^{i\theta})||^2)
\frac{d\theta}{2\pi}\\
& = m_{f^{(1)},\infty}(r) \leq O(\log^+(r\,T_{f,E}(r)))\exc\,\,.
\endsplit
$$
More generally, the estimate
$$
T_{f^{(i)},H^{(i)}}(r) \leq O(\log^+(r\,T_{f,E}(r)))\exc \tag9
$$
holds for each $i\geq 1$. 

The first inequality of Lemma 3.1 implies that
$$
m_{f^{(1)},D^{(1)}}(r) \leq m_{f^{(i)},D^{(i)}}(r)
+O(\log^+(r\,T_{f,E}(r)))\exc \tag10
$$
holds for any $i\geq 1$. The estimate (10) does not give any information on the 
asymptotic behavior of $m_{f^{(1)},D^{(1)}}(r)$, or more generally, of $m_{f^{(i)},D^{(i)}}(r)$. 
However, if we take $i$ sufficiently large, then we can use (10) (combined with (9)) to 
psove an estimate of the form 
$$m_{f^{(i)},D^{(i)}}(r) \leq \ep\,T_{f,E}(r)\exc$$
($\ep$ is any positive number and $E\rightarrow A$ is any fixed ample line bundle). 
We now prove this. 

It is easy to see that $\dim A^{(k)}=2^kn$ and $\dim D^{(k)}=2^k(n-1)$. 
It follows that $D^{(k)}$ has codimension $2^k$. 
There is a natural fibration $D^{(k)} \to D$. We write the fiber over $x\in D$ 
by $D^{(k)}_x$. Fix an arbitrary point $a \in A$. If $2^k>n$, we have
$$
\split
(2^k-1)n & = \text{\rm relative\,\,dimension\,\,of\,\,}A^{(k)} \to A\\
& > 2^k(n-1)=(2^k-1)(n-1)+(n-1)\qquad [\text{\rm since\,\,}2^k>n]\\
& = \dim(D^{(k)})_x+\dim D\\
& = \dim\biggl(
\bigcup_{x\in D}\underbrace{D^{(k)}_x}_{\text{\rm translated\,\,to\,\,a\,\,
fixed\,\,point\,\,}a}\biggr)\,\,.
\endsplit
$$
This implies that if $2^k>n$ there exists a proper subvariety $\wt D^{(k)}$ 
in the fiber of $A^{(k)} \to A$ over $a$ 
such that $D^{(k)} \subset (\wt D^{(k)}\times A)$, 
i.e., $D^{(k)}$ is contained in a ``horizontally flat" proper subvariety. 
Set
$$
\wt D^{(k)}=\bigcup_{x\in D}\underbrace{D_x^{(k)}}_{\text{\rm translated\,\,
to\,\,a\,\,fixed\,\,point\,\,}a}\,\,.
$$
We suppose that $f^{(k)}(\C) \not\subset (\wt D^{(k)}\times A)$. 
Let $E_k$ be a divisor in the fiber of $A^{(k)} \to A$ over $a$ 
chosen so that $\wt D^{(k)}\subset E_{k}$ 
and $f^{(k)}(\C)\not\subset (E_k\times A)$.  
We then have the following estimate:
$$
\split
m_{f^{(k)},D^{(k)}}(r) & \leq m_{f^{(k)},\wt D^{(k)}\times A}(r) 
\leq m_{f^{(k)},E_k\times A}(r)\\
& \leq \exists\,C\,\sum_{i=1}^kT_{f^{(i)},H^{(i)}}(r)\exc\\
& \leq O(\log^+(r\,T_{f,E}(r)))\exc\qquad [\text{\rm from\,\,}(9)]\,\,,
\endsplit
$$
where the constant $C>0$ depends only on the degree of $E_{k}$ in any sense 
in the fiber of $A^{(k)}\to A$, and the degree of $E_{k}$ with the above 
property depends on the degree of $\wt D^{(k)}$. 
For the above argument to make sense, it is sufficient that the 
non-inclusion $f^{(k)}(\C) \not\subset (\wt{D}^{(k)}\times A)$ holds. 
In [K2] I encountered the same situation. The argument given in 
[K2, (2.14), p.147] was not correct in the sense that I overlooked 
this problem (i.e., I ignored constant terms involved in the Lemma on 
logarithmic derivative). 
I would like to take this opportunity to remedy this error. 

Let $\wt D^{(k)}$ be defined as above and we consider a horizontal 
cycle $\wt D^{(k)}\times A$ in $A^{(k)} (\to A)$. 
We assume that the inclusion $f^{(k)}(\C) \subset (\wt D^{(k)}\times A)$ 
holds for all $k$. 
We would like to show that if this happens then the 
holomorphic curve $f:\C \rightarrow A$ is not Zariski-dense. 

We solve this problem by estimating how strongly the jet lift $f^{(k)}:\C\to A^{(k)}$ 
can approximates $D^{(k)}$ for $k$ large (roughly of order $\log_2n$). 
For this purpose, we try to ``deform'' 
the horizontal cycle $\wt D^{(k)}\times A$ 
containing $\wt D^{(k)}\times D$, to non-horizontal ones 
containing $\wt D^{(k)}\times D$ but 
not containing the image $f^{(k)}(\C)$. We will do this in the category 
of $\Q$-cycles. 

We first consider a small non-horizontal ``deformation'' of the horizontal 
divisor $E_k\times A$ (as a $\Q$-divisor). 
We perform such small ``deformation'' by 
introducing a vertical divisor, which means a pull back $F$ of 
some ample divisor of $A$ under the projection $A^{(k)} \to A$, 
and consider a $\Q$-divisor $(E_k\times A)+\ep\,F$ on $A^{(k)}$, 
where $\ep$ is an arbitrary positive rational number (which should be chosen 
to be small)\footnote{\,\,To make the error term $\ep\,T_{f,E}(r)$ in (1) 
smaller, we need to choose $\ep>0$ accordingly smaller in the perturbation, 
although the $\ep$ in perturbation is not exactly the $\ep$ 
in the error term in (1).}. 
Let $d$ be a large positive integer 
such that $L_{k,d,\ep}:=\Cal O_{A^{(k)}}(d\{(E_k\times A)+\ep\,F\})$ 
is a very ample divisor on $A^{(k)}$. 
Then there exists a positive integer $d_0$ such that $d_0/d$ 
is equivalent to a uniform multiple of $\ep$, 
and a holomorphic section $s_{k,d,\ep}$ of $L_{k,d,\ep}$ such that 
the ideal sheaf of the subscheme $(d-d_0)(\wt D^{(k)} \times D)$ 
divides that of $(s_{k,d,\ep})$ (over $D$). 
Thus we get a small ``deformation'' (over $\Q$) $\frac1{d-d_{0}}(s_{k,d,\ep})$ 
of $E_{k}\times A$. 
Choosing an appropriate number of 
such $E_k$'s and taking the intersection of small ``deformation'' (over $\Q$) 
of these $E_k\times A$'s, we get a small non-horizontal ``deformation'' 
of $\wt D^{(k)}\times A$ in the category of $\Q$-cycles. 
Moreover, given ample divisor $D$ (in $A$), $k\in\Bbb Z_{>0}$ and 
any ample line bundle $E$, we have $m(D,k)\in\Bbb Z_{>0}$ 
such that for any $E_{k}$ (in the fiber of $A^{(k)}\to A$) 
containing $\wt D^{(k)}$ and of ``degree 
$\geq m(D,k)$'' w.r.to the ample line bundle $E$, 
we may assume that $\wt D^{(k)}\times A$ and its small ``deformation'' 
(over $\Q$) are disjoint in the fiber of $A^{(k)} \to A$ 
over some point of $A$. 
We now would like to expect the 
non-inclusion $f^{(k)}\not\subset\Supp((s_{k,d,\ep})_0)$. 

However, we cannot expect this non-inclusion because of the following reason. 
Namely, the section $s_{k,d,\ep}\in H^0(A^{(k)},\Cal O(L_{k,d,\ep}))$ 
is special 
in the sense that $(d-d_0)(\wt D^{(k)}\times D)$ 
divides $(s_{k,d,\ep})$ over $D$. So, it might happen 
that the divisor $(s_{k,d,\ep})$ inherits this special property, 
if $\ep$ (the ``strength" of the perturbation) is very small (and we want 
arbitrarily small $\ep$). 
In other words, we have two requirements, i.e., 
\smallskip

(a) to make $\ep$ in 
the perturbation arbitrarily small, and
\smallskip

(b) to find such $\ep$ in the 
perturbation so that the perturbed linear system $|d\{(E_k\times A)+\ep\,F\}|$ 
contains a linear subsystem having exactly $(d-d_0)(\wt D^{(k)}\times D)$ 
as the base scheme. 
\smallskip

However, there may not exist such $\ep>0$ which fulfills both conditions (a) and (b). 

Therefore, the worst such case is that there exists a proper subvariety $Y$ 
of $A$ such that any divisor $(s_{k,d,\ep})$ having the above divisibility 
property necessarily contains a subscheme of the form $Z_{k} \times Y$ 
($Z_{k}$ being a subvariety of $\wt D^{(k)}$) 
and the inclusion $f^{(k)}(\C) \subset Z_{k} \times Y$ holds. 
We cannot avoid this possibility in the small ``deformation'' 
of $E_k\times A$ (see the discussion after Theorem 3.1 below). 
Therefore, instead of looking for a ``good'' $E_k \times A$, 
we now put the assumption that  the assumption 
the holomorphic curve $f:\C \rightarrow  A$ is 
{\bf Zariski-dense (algebraically non-degenerate)}.

Under this assumption, we have the estimate
$$
m_{f^{(k)},(d-d_0)\,D^{(k)}}(r) \leq m_{f^{(k)},(s_{k,d,\ep})}(r)\,\,.
$$
Dividing the both sides by $d-d_0$ we have
$$
m_{f^{(k)},D^{(k)}}(r) 
\leq \frac{d}{d-d_0} T_{f^{(k)},(E_k\times A)+\ep\,F}(r)\,\,.
$$
If we choose $\ep$ sufficiently small and $d$ sufficiently large, we can 
make $d_0/d$ arbitrarily small. Combining this with (9), we finally have 
the estimate of type
$$
m_{f^{(k)},D^{(k)}}(r) \leq \ep\,T_{f,E}(r)\exc\,\,,
$$
where $E$ is any fixed ample line bundle on $A$ and $\ep$ is any positive 
number. We have thus proved the main result in Yamanoi's paper [Y2] (the 
Second Main Theorem for Zariski-dense holomorphic curves into Abelian varieties) 
by studying the truncated counting functions 
in the framework of the Radon transformation. 

\proclaim{Theorem 3.1 (Second Main Theorem with truncation level 1)} 
Let $A$ be an Abelian variety, $D$ any reduced divisor 
and $\ep$ any small positive number. 
Let $f:\C \to A$ be a Zariski-dense holomorphic curve. 
Then we have
$$
m_{f,D}(r)+N^1_{f,D}(r) \leq \ep\,T_{f,E}(r)\exc\,\,.
$$
\endproclaim
\noindent

The assumption $f:\C \rightarrow A$ being Zariski dense cannot be removed. 
Indeed, we choose an Abelian variety $A$ which admits a non-trivial proper 
Abelian subvariety $i:B\hookrightarrow A$ and a reduced divisor $D$ of $A$ 
having the property that $i^{*}D=m\,D_{B}$ for $\Z \ni m \geq 2$ and $D_{B}$ 
is a reduced divisor on $B$. Let $f:\C \rightarrow A$ be a holomorphic curve 
such that $f(\C) \subset B$. This choice of $A$ and $f$ violates the assumption 
of Theorem 3.1. 
On the other hand, since all intersections 
of $f$ and $D$ are with multiplicity $\geq 2$, 
we have $N^{1}_{f,D}(r)\geq \frac12 N_{f,D}(r)$. This implies
$$
\split
m_{f,D}(r)+N^{1}_{f,D}(r) & 
\geq \frac12 m_{f,D}(r)+\frac12(m_{f,D}(r)+N_{f,D}(r))\\
& = \frac12 m_{f,D}(r)+\frac12 T_{f,D}(r)+O(1)\\
& \geq \frac12 T_{f,D}(r)\,\,.
\endsplit
$$
This shows that the Second Main Theorem with the level 1 residual counting 
function does not hold for this choice of $A$ and $f$. Therefore 
the assumption $f$ being Zariski dense cannot be removed. 
\medskip

To sum up, we have proved the following: if $f:\C \to A$ is Zariski-dense, then 
we conclude that the intersection points of $f$ and $D$ with high multiplicity are ``rare" 
in the sense that Theorem 3.1 holds.

\enddocument